\documentclass[11pt]{article}
\usepackage[T2A]{fontenc}

\usepackage[russian,english]{babel}
\usepackage[russian]{babel}
\usepackage{amscd,amssymb,amsmath,amsthm}
\usepackage[arrow,matrix]{xy}
\usepackage{graphicx}
\usepackage{color}
\usepackage{cite}

\theoremstyle{plain}
\newtheorem{theorem}{Теорема}

\newtheorem{definition}{Определение}

\newtheorem{remark}{Замечание}

\begin{document}

\sloppy
\begin{center}
{\bf \large {О НОВЫХ СЛАБО ПЕРИОДИЧЕСКИХ ГИББСОВСКИХ МЕРАХ МОДЕЛИ
ИЗИНГА НА ДЕРЕВЕ КЭЛИ}}
\end{center}

\begin{center}
 М.М.Рахматуллаев \footnote{Институт математики, ул. Дурмон йули, 29, Ташкент, 100125, Узбекистан.\\
E-mail: mrahmatullaev@rambler.ru}
\end{center}
\selectlanguage{russian}

\begin {abstract} Рассматривается модель Изинга на дереве Кэли. Для модели
Изинга найдены новые слабо периодические меры, соответствующие
нормальным делителям индекса 2 группового представления дерева
Кэли.
\end {abstract}

\textbf{Ключевые слова:} дерево Кэли, мера Гиббса, модель Изинга,
слабо периодические меры.

УДК.517.98+530.1

\selectlanguage{english}
\begin {abstract}We consider Ising model on the Cayley tree. For the Ising model, we find a new weakly periodic Gibbs
measures corresponding to normal subgroups of indices two in the
group representation of a Cayley tree.\end {abstract}


\selectlanguage{russian}

\begin{center}
\textbf{1. Введение}
\end{center}

Каждой мере Гиббса сопоставляется одна фаза физической системы.
Если существует более чем одна мера Гиббса, то говорят, что
существуют фазовые переходы. Основная проблема для данного
гамильтониана - это описание всех отвечающих ему предельных мер
Гиббса. Эти меры, в основном, были трансляционно-инвариантными,
либо периодическими с периодом два. Более того, для многих моделей
на дереве Кэли доказано, что множество периодических мер Гиббса
очень бедно, т.е. существуют периодические гиббсовские меры с
периодом два. В работах \cite{1}-\cite{3} для модели Изинга
описаны трансляционно-инвариантные меры Гиббса на дереве Кэли.
Описанию периодических гиббсовских мер для некоторых моделей с
конечным числом радиуса взаимодействия, которые, в основном, были
либо трансляционно-инвариантными, либо периодическими с периодом
два, посвящены работы (\cite{4},\cite{5}), .

Чтобы получить более широкое множество гиббсовских мер, в работах
\cite{7}-\cite{10} введены более общие понятия периодической меры
Гиббса, т.е. слабо периодические гиббсовские меры и доказано
существование таких мер для модели Изинга на дереве Кэли. В
работах \cite{1}, \cite{6}, \cite{10} изучены континуальные
множества непериодических мер Гиббса для модели Изинга на дереве
Кэли.

В работах \cite{7} и \cite{8} на некоторых инвариантных множествах
при некоторых условиях на параметры найдены слабо-периодические
(не периодические) меры Гиббса для модели Изинга на дереве Кэли.
Настоящая работа является продолжением работ \cite{7} и \cite{8}.
Доказывается существование новых мер Гиббса, отличных от мер,
описанных в работах \cite{7} и \cite{8}.

Структура работы: В пункте 2 даются основные определения и
постановка задачи. Пункт 3 посвящен слабо периодическим мерам
Гиббса.

\begin{center}
\textbf{2. Основные определения и постановка задачи}
\end{center}

Пусть $\tau^{k}=(V,L), k\geq1$ есть дерево Кэли порядка $k$, т.е.
бесконечное дерево, из каждой вершины которого выходит равно $k+1$
ребер, где $V$ - множество вершин, $L$ - множество ребер $\tau^k.$

Пусть $G_k$ - свободное произведение
 $k+1$ циклических групп $\{e,a_i\}$ второго порядка с образующими $a_1, a_2, ...
 a_{k+1}$, соответственно т.е. $a_i^2=e$.

Существует взаимно-однозначное соответствие между множеством
вершин $V$ дерево Кэли порядка $k$ и группой $G_k$(см. \cite{1},
\cite{11}).

Это соответствие строится следующим образом. Произвольной
фиксированной вершине $x_0\in V$ поставим соответствие единичный
элемент $e$ группы $G_k$. Так как рассматриваемый граф без
ограничения общности можно считать плоским, то каждой соседной
вершине точки $x_0$ (т.е. $e$) поставим в соответствие образующую
$a_i, i=1,2,...k+1$ по положительному направлению (см. рис.1).

\begin{figure}
  \includegraphics[width=0.8\textwidth]{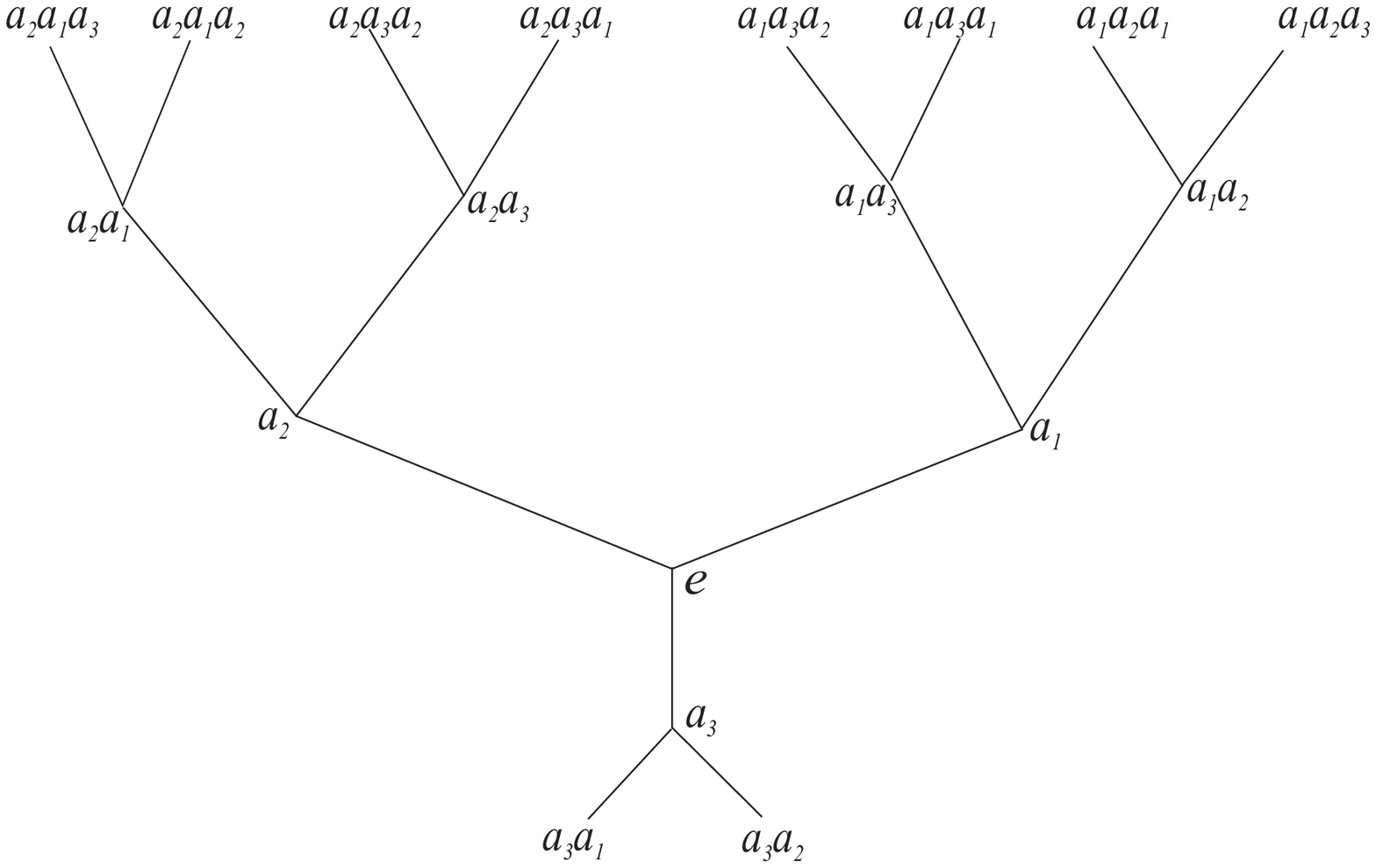}\\
  \caption{\textit{Дерево Кэли $\tau^2$ и элементы группового представления вершин.}}\label{1}
\end{figure}

Теперь в каждой вершине $a_i$ определим слово длины два $a_ia_j$
соседних вершин $a_i$. Так как одна из соседних вершин вершины
$a_i$ есть $e$, то положим $a_ia_i=e$ и тогда нумерация остальных
соседних вершин $a_i$ производится однозначно по вышеприведенному
правилу нумерации.  Далее, для соседних вершин вершины $a_ia_j$
определим слово длины три следующим образом. Так как одна из
соседних для $a_ia_j$  вершин есть $a_i$, то положим
$a_ia_ja_j=a_i$ и тогда нумерация остальных соседних вершин
производится однозначно и имеют вид $a_ia_ja_l, i,j,l=1,2,...k+1.$
Это соответствие согласуется с предыдущим шагом, так как
$a_{i}a_ja_j=a_ia_j^2=a_i$. Таким образом, можно установить
взаимно-однозначное соответствие между множеством вершин дерева
Кэли $\tau^k$ и группой $G_k$.

Представление, построенное выше, называется правым, т.к. в этом
случае если $x$ и $y$ соседние вершины, а $g$ и $h\in G_k$
соответствующие им элементы группы, тогда либо $g=ha_i$, либо
$h=ga_j$ для некоторых $i$ или $j$. Аналогично определяется левое
представление.

Рассмотрим в группе $G_k$ (соответственно на дереве Кэли)
преобразование левого (правого) сдвига, определяемые следующим
образом: для $g\in G_k$ положим
$$
T_{g}(h)=gh,  (T_{g}(h)=hg), \quad \forall h\in G_k.
$$
Совокупность всех левых (правых) сдвигов на $G_k$ изоморфна группе
$G_k$.

Любое преобразование $S$ группы $G_k$ индуцирует преоброзование
$\widehat{S}$ на множестве вершин $V$ дерева Кэли $\tau^k$.
Поэтому мы отоджествляем $V$ и $G_k$.

\begin{theorem}
Группа левых (правых) сдвигов на правом (левом) представлении
дерева Кэли является группой трансляций.(см. \cite{1}, \cite{11})
\end{theorem}
Для произвольной точки $x^0\in V$ положим $ W_n=\{x\in V|
d(x^0,x)=n\}$, $V_n=\bigcup\limits_{m=0}^nW_m,$ $L_n=\{<x,y>\in L|
x,y\in V_n\}$, где $d(x,y)$-расстояние между  $x$ и $y$ на дереве
Кэли, т.е. число ребер пути, соединяющее  $x$ и $y$ .

Пусть $\Phi=\{-1,1\}$ и $\sigma\in \Omega=\Phi^V$ конфигурация, то
есть $\sigma=\{\sigma(x)\in \Phi: x\in V\}$. Пусть $A\subset V$.
Обозначим через $\Omega_A$ пространство конфигураций, определенных
на множестве $A$ и принимающих значения из $\Phi=\{-1,1\}$.

Рассмотрим гамильтониан модели Изинга

$$
H(\sigma)=-J\sum _{<x,y>\in L}\sigma(x) \sigma(y), \eqno (1)
$$
где $J\in R, <x,y>$ - ближайшие соседи. Заметим, что элементы
$x,y$ группы $G_k$ являются соседними, если существует образующая
$a_i$ т.е. $x=ya_i$.

Пусть $h_x\in R, x\in V.$ Для каждого $n$ определим меру $\mu_n$
на  $\Omega_{V_n}$, полагая
$$
\mu_{n}(\sigma_n)=Z^{-1}_n \exp\{-\beta H(\sigma_n)+\sum_{x\in
W_n}h_x \sigma (x)\},\eqno(2)
$$ где
$\beta=\frac{1}{T}$ ($T>0$-температура),
$\sigma_{n}=\{\sigma(x),x\in V_n\}\in \Omega_{V_n}$, $Z^{-1}_n -$
нормирующий множитель и $$H(\sigma_n)=-J\sum_{<x,y>\in L_n}\sigma
(x) \sigma (y).$$ Условие согласованности для $\mu_n(\sigma_n),
n\geq 1,$ определяется равенством
$$
\sum_{\sigma^{(n)}}\mu_n(\sigma_{n-1},\sigma^{(n)})=\mu_{n-1}(\sigma_{n-1}),\eqno(3)
$$
где $\sigma^{(n)}=\{\sigma(x),x\in W_n\}.$

Пусть $\mu_n, n\geq1$ - последовательность мер на $\Omega_{V_n}$,
обладающая свойством согласованности (3). Тогда в силу теоремы
Колмогорова существует и притом единственная предельная мера $\mu$
на $\Omega_{V}=\Omega$ (которая называется предельной мерой
Гиббса) такая, что для каждого $n=1,2,...$
$$\mu(\sigma_n)=\mu_n(\sigma_n).$$
 Известно, что меры (2) удовлетворяют (3) тогда и только тогда,
 когда совокупность величин $h=\{h_x,
x\in G_k \}$ удовлетворяет
$$h_x=\sum_{y \in S(x)}f(h_y,\theta), \eqno (4)$$ где $S(x)$-
множество  "прямых потомков" точки $x\in V$ и
$f(x,\theta)=arcth(\theta thx)$,  \\ $\theta=th(J\beta)$,
 $\beta=\frac{1}{T}, T>0$ температура (см.\cite{1},
\cite{2},\cite{3}).\

\begin{definition} Совокупность величин $h=\{h_x,x\in G_k\}$ называется $
\widehat{G}_k$-\textit{периодической}, если $h_{xy}=h_x$ для
$\forall x\in G_k, y\in\widehat{G}_k$. $G_k$-периодическая мера
называется \textit{трансляционно-инвариантной}.
\end{definition}

Для $ x\in G_k $ обозначим через $x_{\downarrow} =\{y\in
G_k:<x,y>\}\backslash S(x)$.

Пусть $G_k/\widehat{G}_k=\{H_1,...,H_r\}$ - фактор группа, где
$\widehat{G}_k$ - нормальный делитель индекса $r\geq 1.$

\begin{definition} Совокупность величин $h=\{h_x,x\in G_k\}$
называется $\widehat{G}_k$ - \textit{слабо периодической}, если
$h_x=h_{ij}$ при $x\in H_i, x_{\downarrow}\in H_j$ для $\forall
x\in G_k$.
\end{definition}

\begin{definition} Мера $\mu$ называется
$\widehat{G}_k$-\textit{(слабо) периодической}, если она
соответствует $\widehat{G}_k$-(слабо) периодической совокупности
величин $h$.
\end{definition}

Целью работы является описание новых слабо периодических
гиббсовских мер для модели Изинга, отличных от мер, которые
получены в \cite{7} и \cite{8}.

\vskip 0.4 truecm

\begin{center}
\textbf{3. Слабо периодические меры}
\end{center}

Заметим, что любой нормальный делитель индекса 2 группы $G_k$
имеет вид $H_A=\{x\in G_k:\sum\limits_{i\in A}w_x(a_i)-$ четно$
\}$, где $ \emptyset \neq A\subseteq N_k=\{1,2,...,k+1\}$ и
$w_x(a_i)-$ число букв $a_i$ в слове $x\in G_k$ (см.[4]). Опишем
$H_A$-слабо периодические меры Гиббса. Заметим, что в случае
$|A|=k+1$, т.е. $A=N_k$ понятие слабо периодичности совпадает с
обычной периодичностью, где $|A|$ обозначает число элементов
множества $ A $. Поэтому рассмотрим $A\subset N_k$ такое, что
$A\neq N_k$. Тогда в силу (4) $H_A$-слабо периодическая
совокупность $h$ имеет вид

$$
h_x=\left\{%
\begin{array}{ll}
    h_{1}, & {x \in H_A, \ x_{\downarrow} \in H_A} \\
    h_{2}, & {x \in H_A, \ x_{\downarrow} \in G_k \backslash H_A} \\
    h_{3}, & {x \in G_k \backslash H_A, \ x_{\downarrow} \in H_A} \\
    h_{4}, & { x \in G_k \backslash H_A, x_{\downarrow}  \in G_k \backslash H_A,}    \\
\end{array}%
\right.\eqno(5)
$$
где $h_{i},i=\overline{1,4}$, удовлетворяет системе уравнений:
$$
\left\{%
\begin{array}{ll}
    h_{1}=|A|f(h_{3},\theta)+(k-|A|)f(h_{1},\theta) \\
    h_{2}=(|A|-1)f(h_{3},\theta)+(k+1-|A|)f(h_{1},\theta) \\
    h_{3}=(|A|-1)f(h_{2},\theta)+(k+1-|A|)f(h_{4},\theta) \\
    h_{4}=|A|f(h_{2},\theta)+(k-|A|)f(h_{4},\theta).\\
\end{array}%
\right.\eqno(6)
$$

Рассмотрим отображение $W:R^4\rightarrow R^4$, определенное
следующим образом:
$$
\left\{%
\begin{array}{ll}
    h'_{1}=|A|f(h_{3},\theta)+(k-|A|)f(h_{1},\theta) \\
    h'_{2}=(|A|-1)f(h_{3},\theta)+(k+1-|A|)f(h_{1},\theta) \\
    h'_{3}=(|A|-1)f(h_{2},\theta)+(k+1-|A|)f(h_{4},\theta) \\
    h'_{4}=|A|f(h_{2},\theta)+(k-|A|)f(h_{4},\theta).\\
\end{array}%
\right.\eqno(7)
$$
\\
Заметим, что (6) есть уравнение $h=W(h)$.

Легко показать, что отображение $W$ имеет инвариантные множества
следующих видов:

$$ I_1 =\{h\in R^4: h_1=h_2=h_3=h_4\}, \ \ I_2 =\{h\in R^4:
h_1=h_4; h_2=h_3\}$$ $$ I_3 =\{h\in R^4: h_1=-h_4; h_2=-h_3\}. $$

Обозначим $\alpha=\frac{1-\theta}{1+\theta}.$ В работе \cite{7}
доказана следуюшая

\begin{theorem}  1) Для  модели Изинга все
$H_A$-слабо периодические меры Гиббса на $I_1, I_2$ являются
трансляционно-инвариантными.

2) При $|A|=k$ и $\theta >0$ все $H_A$ - слабо периодические меры
Гиббса являются трансляционно-инвариантными.

3) При $|A|=1,k=4$ существует критическое значение
$\alpha_{cr}(\approx 0,152)$ такое, что при $0<\alpha<\alpha_{cr}$
существуют пять $H_A$-слабо периодические меры Гиббса $\mu_0,
\mu^-_1, \mu^+_1, \mu^-_2, \mu^+_2$; при $\alpha=\alpha_{cr}$
существуют три $H_A$-слабо периодические меры Гиббса $\mu_0,
\mu^-_1, \mu^+_1;$ при $\alpha>\alpha_{cr}$ существует
единственная $H_A$-слабо периодическая мера Гиббса $\mu_0.$

4) При $|A|=1, k>5$ и $\theta\in (\theta_1 , \theta_2)$ на $I_3$
существуют три $H_A$ - слабо периодические меры Гиббса
$\mu^-,\mu^0,\mu^+$ где
$\theta_{1,2}=\frac{k-1\pm\sqrt{k^2-6k+1}}{2k}$.
\end {theorem}

Во втором пункте теоремы 2 рассматривается случай $\theta >0$.
Поэтому естественно рассмотреть случай $\theta <0 (\alpha>1)$.

Пользуясь тем, что $f(h,\theta)=arcth(\theta th
h)=\frac{1}{2}\ln\frac{(1+\theta)e^{2h}+(1-\theta)}{(1-\theta)e^{2h}+(1+\theta)}$
и обозначая $z_i=e^{2h_i}$, где $ (i=\overline{1,4})$, из (7)
получим следующую систему уравнений
$$
\left\{%
\begin{array}{ll}
    z_{1}=(\frac{z_3+\alpha}{\alpha z_3+1})^{|A|}\cdot(\frac{z_1+\alpha}{\alpha z_1+1})^{(k-|A|)} \\
    z_{2}=(\frac{z_3+\alpha}{\alpha z_3+1})^{|A|-1}\cdot(\frac{z_1+\alpha}{\alpha z_1+1})^{(k+1-|A|)}\\
    z_{3}=(\frac{z_2+\alpha}{\alpha z_2+1})^{|A|-1}\cdot(\frac{z_4+\alpha}{\alpha z_4+1})^{(k+1-|A|)} \\
    z_{4}=(\frac{z_2+\alpha}{\alpha z_2+1})^{|A|}\cdot(\frac{z_4+\alpha}{\alpha z_4+1})^{(k-|A|)}.\\
\end{array}%
\right.\eqno(8)
$$

Введем следующее обозначение $f(x)=\frac{x+\alpha}{\alpha x+1}$.
Тогда система уравнений (8) при $|A|=k$ сводится к следующей
системе уравнений
$$
\left\{%
\begin{array}{ll}
    z_{1}=\left(f(z_3)\right)^{k}\\
    z_{2}=\left(f(z_3)\right)^{k-1}\cdot \left(f(z_1)\right)\\
    z_{3}=\left(f(z_2)\right)^{k-1}\cdot \left(f(z_4)\right)\\
    z_{4}=\left(f(z_2)\right)^{k}.
\end{array}%
\right.\eqno(9)$$

Справедлива следующая

\begin{theorem} Пусть $|A|=k,\alpha>1.$ \\
1) При $k\leq 3$ все $H_A-$ слабо периодические меры Гиббса на $I_3$ являются трансляционно-инвариантными. \\
2) При $k=4$ существует критическое значение $\alpha_{cr}(\approx
6,3716)$, такое, что при $\alpha<\alpha_{cr}$  на $I_3$ существует
одна $H_A-$ слабо периодическая мера Гиббса; при
$\alpha=\alpha_{cr}$ на $I_3$ существует три $H_A-$ слабо
периодические меры Гиббса; при $\alpha>\alpha_{cr}$  на $I_3$
существует пять $H_A-$ слабо периодические меры Гиббса.
\end{theorem}
\begin{proof} Система уравнений (9) на инвариантном множестве
$I_3$ имеет следующий вид:
$$\left\{%
\begin{array}{ll}
    z_{1}=\left(f(\frac{1}{z_2})\right)^{k}\\
    z_{2}=\left(f(\frac{1}{z_2})\right)^{k-1}\cdot \left(f(z_1)\right).
\end{array}%
\right.\eqno(10)
$$
Подставляя первое уравнение (10) во второе, получим
$$
z_2=\left(\frac{1+\alpha
z_2}{\alpha+z_2}\right)\frac{\alpha(\alpha+z_2)^k+(1+\alpha
z_2)^k}{(\alpha+z_2)^k+\alpha(1+\alpha z_2)^k}.\eqno(11)
$$
Введем обозначение $u=\frac{z_2+\alpha}{\alpha z_2+1}$. Тогда
уравнение (11) сводится к уравнению
$$
u^{2k}-\alpha u^{2k-1}+\alpha^2 u^{k+1}-\alpha^2 u^{k-1}+\alpha
u-1=0. \eqno(12)
$$
Заметим, что уравнение (12) можно переписать в виде
$$
(u^2-1)P_{2k-2}(u)=0, \eqno(13)
$$
где $P_{2k-2}(u)-$ симметрический многочлен степени $2k-2$.
Обозначая $u+\frac{1}{u}=\xi$, степень уравнения $P_{2k-2}(u)=0$
можно уменьшить в два раза, т.е. привести к виду $P_{k-1}(\xi)=0$.
Но при $k\geq 6$ уравнение не решается в кубатуре.

Сначало докажем второй пункт теоремы 2. Пусть $k=4.$ В этом случае
уравнение (13) имеет следующий вид
$$
(u^2-1)\left(u^6-\alpha u^5+u^4+(\alpha^2-\alpha)u^3+u^2-\alpha
u+1\right)=0. \eqno(14)
$$
Отсюда $ u^2-1=0$ или
$$
u^6-\alpha u^5+u^4+(\alpha^2-\alpha)u^3+u^2-\alpha u+1=0.
\eqno(15)
$$
Так как $u>0$, то $u=1$ является одним из решений уравнения (14)
для любого $\alpha>0$. Рассмотрим уравнение (15). Обозначив $\xi
=u+\frac{1}{u}>2$, получим следующее уравнение
$$
\xi^3-\alpha \xi^2-2\xi+\alpha^2+\alpha=0.
$$
Рассмотрим функцию $\gamma(\xi)=\xi^3-\alpha
\xi^2-2\xi+\alpha^2+\alpha=0$ и изучим ее. Очевидно, что
$\gamma(2)>0$ при всех $\alpha>1$. Из $\gamma'(\xi)=0$ получим
$\xi=\xi_*=\frac{\alpha+\sqrt{\alpha^2+6}}{3}$. Из $\xi_*>2$
находим, что $\alpha>2,5$. Если $\gamma(\xi_*)<0$, то уравнение
$\gamma(\xi)=0$ имеет точно два решения $\xi_1>2, \xi_2>2$. Из
$u+\frac{1}{u}=\xi$ получим четыре решения уравнения (15). Таким
образом, в случае $\alpha>\alpha_{cr}$ уравнение (14) имеет 5
решений, где критическое значение $\alpha_{cr}(\approx 6,3716)$
находится из уравнения $\gamma(\xi_*)=0$ т.е.
$$
9\alpha+27\alpha^2-2\alpha^3=2(\sqrt{\alpha^2+6})^3.
$$
Рассмотрим
$\psi(\alpha)=2\alpha^3-27\alpha^2-9\alpha+2(\sqrt{\alpha^2+6})^3,
\alpha>2,5.$ Заметим, что $\psi(\alpha)$ возрастает в интервале
$(\alpha', +\infty)$, где $\alpha'=\frac{48+\sqrt{264}}{12}>5.$
Так как $\psi(5)<0$, из монотонности $\psi$ следует существование
и единственность $\alpha_{cr}$ такое, что $\psi(\alpha_{cr})=0.$ В
случае $\alpha=\alpha_{cr}$ уравнение $\gamma(\xi)=0$ имеет
единственное решение $\xi=\xi_*>2$. Из $u+\frac{1}{u}=\xi$ получим
два решения уравнения (15). Следовательно уравнение (14) имеет три
решения. А в случае $\alpha<\alpha_{cr}$ уравнение $\gamma(\xi)=0$
не имеет решений $\xi=\xi_*>2$. Следовательно, в этом случае
уравнение (14) имеет единственное решение $u=1$.

Теперь докажем первый пункт теоремы 2. Пусть $k=3$. Тогда
уравнение (13) имеет следующий вид
$$
(u^2-1)(u^4-\alpha u^3+(\alpha^2+1)u^2-\alpha u+1)=0. \eqno(16)
$$

В этом случае тоже $u=1$ является одним из решений уравнения (16).
Поэтому предположим, что $u\neq1.$ Тогда из (16) имеем
$$
u^4-\alpha u^3+(\alpha^2+1)u^2-\alpha u+1=0.
$$
Обозначая $\xi=u+\frac{1}{u}$, получим следующее квадратное
уравнение
$$
\xi^2+\alpha \xi +\alpha^2-1=0.
$$
Легко показать, что это уравнение при $\xi>2$ и $\alpha>1$ не
имеет решения, следовательно, уравнение (16) не имеет решения,
кроме $u=1$.

Аналогично рассуждая, легко доказать, что при $k\leq 2$ уравнение
(12) не имеет решения, кроме $u=1.$ Теорема доказана.\end{proof}

\begin{remark} Заметим, что одна из мер, найденных в пункте 2 теоремы 3,
является трансляционно-инвариантной, а все остальные являются
$H_A-$ слабо периодическими (не периодическими), которые
отличаются от мер, найденных в \cite{7} и \cite{8}.
\end{remark}

\begin{remark} Компьютерный анализ показывает, что
утверждения пункта 2 теоремы 3 верны для всех $k\geq 4$, т.е.
независимо от значения $k$ на инвариантном множестве $I_3$
существуют не более пяти $H_A-$ слабо периодических мер
Гиббса.\end{remark}

\begin{remark} В \cite{8} при одном нормальном делителе индекса
4 получены не менее 7 слабо периодических мер Гиббса, где
критическое значение $\alpha$ совпадает с
$\alpha_{cr}$.\end{remark}

\begin{figure}
  \includegraphics[width=0.5\textwidth]{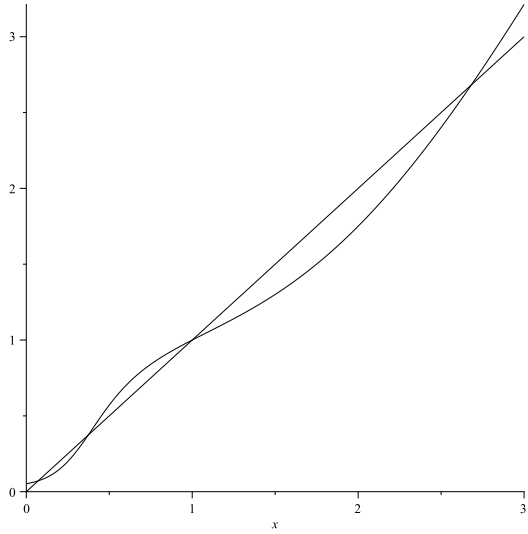}\\
  \caption{\textit{График функции $\varphi(x)$ и $y=x $ при $k=5 \ \ \alpha=3.$ }}\label{1}
\end{figure}

\begin{theorem} Пусть $|A|=k, \ \  k\geq 6$, $\alpha>1 $
 и $\alpha_{1,2}=\frac{k-1\pm\sqrt{k^2-6k+1}}{2}$. Тогда верны
 следуюшие утверждения: \\
1) При $\alpha\in (\alpha_1 , \alpha_2)$ существуют не менее трех
$H_A$ - слабо периодических мер Гиббса.\\
2) При $\alpha\notin (\alpha_1 , \alpha_2)$ существует не менее
одной $H_A$ - слабо периодической меры Гиббса.
\end{theorem}

\begin{proof}  Первое и четвертое уравнения из (9) подставим в
третье и во второе уранения той же системы, соответственно. А в
полученной системе уравнений третье уравнение подставим во второе
и получим следуюшую систему уравнений:
$$
\left\{%
\begin{array}{ll}
    z_{1}=\left(f(z_3)\right)^{k}\\
    z_{2}=\left(f\left(\left(f(z_2)\right)^{k-1}\cdot f\left(\left(f(z_2)\right)^{k}\right)\right)\right)^{k-1}\cdot f\left(\left(f\left(\left(f(z_2)\right)^{k-1}\cdot f\left(\left(f(z_2)\right)^{k}\right)\right)\right)^{k}\right)\\
    z_{3}=\left(f(z_2)\right)^{k-1}\cdot f(\left(f(z_2)\right)^{k})\\
    z_{4}=\left(f(z_2)\right)^{k}.
\end{array}%
\right.$$

Введем следующее обозначение:
$$\varphi(x)=\left(f\left(\left(f(x)\right)^{k-1}\cdot f\left(\left(f(x)\right)^{k}\right)\right)\right)^{k-1}\cdot
f\left(\left(f\left(\left(f(x)\right)^{k-1}\cdot
f\left(\left(f(x)\right)^{k}\right)\right)\right)^{k}\right).\eqno
(17)$$

Известно, что функция $f(x)$ ограничена, отсюда следует, что и
$\varphi (x)$ тоже будет ограниченной. Ясно, что
$$
\varphi(0)>0;\quad \varphi(1)=1.
$$
Рассмотрим производную функции $\varphi(x)$:
$$
\varphi'(x)=  { \left( f \left(  \left( f \left( x \right) \right)
^{k-1}f
 \left(  \left( f \left( x \right)  \right) ^{k} \right)  \right)
 \right) ^{k-2} f' \left(  \left( f \left( x
 \right)  \right) ^{k-1}f \left(  \left( f \left( x \right)  \right) ^
{k} \right)  \right)  \left( f \left( x \right)  \right) ^{k-2}
  f'(x)}\times
$$$$ \Bigg[ kf \left(
 \left( f \left( x \right)  \right) ^{k} \right) -f \left(  \left( f
 \left( x \right)  \right) ^{k} \right) + kf'\left(  \left( f \left( x \right)  \right) ^{k} \right)  \left( f
 \left( x \right)  \right) ^{k} \Bigg] \times$$
 $$  \Bigg\{ kf \left(  \left( f
 \left(  \left( f \left( x \right)  \right) ^{k-1}f \left(  \left( f
 \left( x \right)  \right) ^{k} \right)  \right)  \right) ^{k}
 \right) -f \left(  \left( f \left(  \left( f \left( x \right)
 \right) ^{k-1}f \left(  \left( f \left( x \right)  \right) ^{k}
 \right)  \right)  \right) ^{k} \right) +$$
 $$ k f'\left(  \left( f \left(  \left( f \left( x \right)  \right) ^{k-1}f
 \left(  \left( f \left( x \right)  \right) ^{k} \right)  \right)
 \right) ^{k} \right)  \left( f \left(  \left( f \left( x \right)
 \right) ^{k-1}f \left(  \left( f \left( x \right)  \right) ^{k}
 \right)  \right)  \right) ^{k}\Bigg\}. $$

Учитывая $$f(1)=1, f'(1)={\frac {1-\alpha}{1+\alpha}},$$ имеем \\
$$ \varphi'(1)=\frac { \left( \alpha-1 \right) ^{2} \left( \alpha+1-2k \right) ^{
2}}{ \left( 1+\alpha \right) ^{4}}.\eqno(18)$$

При $ \varphi'(1)>1$ графики функций $ y=\varphi(x)$ и $y=x$ как
минимум пересекаются трижды, т.е. уравнение (17) имеет не менее
трех решений. Поэтому будем решать следующее неравенство:
$$
\varphi'(1)>1. \eqno(19)
$$
Из (18) и (19) получим

$$
\frac{4(\alpha-\frac{k-1}{k+1})(\alpha-\alpha_1)(\alpha-\alpha_2)
}{(k+1) \left( 1+\alpha \right) ^{4}}<0,
$$
где $\alpha_1, \alpha_2$ определены в теореме 3.

Ясно, что если $k<5$, то $\alpha_1, \alpha_2$ являются
комплексными, а если $k=5$, то $\alpha_1=\alpha_2$ и неравенство
(19) имеет решение
$$(-\infty, \frac{k-1}{k+1}),$$ но
$\alpha$ должна принадлежать интервалу $(1, +\infty)$.

Теперь рассмотрим $k\geq 6,$ тогда неравенство (19) имеет решение
$( \alpha_1, \alpha_2)$. В этом случае система уравнений (16)
имеет три решения, т.е. существует три $H_A-$слабо периодические
меры Гиббса. При $\alpha\notin( \alpha_1, \alpha_2)$ прямая $y=x$
и график функции $y=\varphi (x)$ пересекаются в точке $(1;1)$ т.е.
существует не менее одной $H_A-$ слабо-периодической меры Гиббса.
Теорема доказана. \end{proof}

\begin{remark} Заметим, что одна из мер, найденных в
теореме 4, является трансляционно-инвариантной, а все остальные
являются $H_A-$ слабо периодическими (не периодическими), которые
являются новыми.\end{remark}

\begin{remark} Также заметим, что при некоторых
$\alpha\in(1,+\infty)$ появляются новые $H_A-$ слабо периодические
меры Гиббса, отличные от мер из Теоремы 4. Например, при $k=5, \ \
\alpha=3$ из рисунка 2 и из ограниченности функции $ y=\varphi(x)$
можно убедиться, что существуют пять $H_A-$ слабо периодических
мер Гиббса, одна из которых является
трансляционно-инвариантной.\end{remark}

\begin{remark} Новые меры Гиббса, описанные в теоремах 3 и
4, дают возможность описать континуум множества непериодических
гиббсовских мер, отличных от ранее известных.\end{remark}

\vskip 0.4 truecm

\begin{center}
\textbf{4. Обсуждение}
\end{center}

Широко распространена (см. \cite{*1}-\cite{*4}) формулировка
модели Изинга, согласно которой каждому узлу $x$ решетки
сопоставляется переменная $\sigma(x)$, принимающая численные
значения $+1$ или $-1$.

Если "объекты", связанные с узлами $x$ и $x'$, находятся в одном
состоянии, то $\sigma(x)\sigma(x')=+1$, а если в разных, то
$\sigma(x)\sigma(x')=-1$. Ясно, что в такой интерпретации понятие
"объекта", связанного с узлом $x$, может трактоваться весьма
широко. Это могут быть, например, магнитный момент иона в
кристалле, имеющей два возможных направления \cite{*15} или атомы
двух сортов в бинарном сплаве \cite{*16}(значение $\sigma(x)=+1$
соответсвует занятию $x$-го узла атомом одного сорта, а
$\sigma(x)=-1$-атомом другого сорта). Другие интерпретации модели
Изинга связаны исследованием явления адсорбции на поверхности
\cite{*17}, плавлением ДНК \cite{*18}, теорией решеточного газа
\cite{*19} и целым рядом других вопросов теории фазового перехода
типа порядок-безпорядок.

Для рассматриваемой модели на $\tau^4$ найдено критическое
значение $\alpha_{cr}(\approx 6,3716)$, такое, что при
$\alpha=\alpha_{cr}$ существует три $H_A-$ слабо периодические
меры Гиббса; при $\alpha>\alpha_{cr}$ существуют пять $H_A-$ слабо
периодические меры Гиббса. Также для этой модели на дереве Кэли
$\tau^k$ при $k\geq 6$ доказано существование не менее трех слабо
периодических мер Гиббса. Различные гиббсовские меры соответствуют
различным фазам системы. В этом случае существуют фазовые
переходы.

Заметим, что построенные в наших случаях фазы отличаются от
известных фаз этой модели, что они соответствуют слабо
периодическим мерам Гиббса (cм.\cite{GRS}, \cite{7}-\cite{10}).

\textbf{Благодарность.} Автор выражает глубокую признательность
профессору У. А. Розикову за постановку задачи и полезные советы
по работе.

\newpage
\begin{center}
\textbf{Сведения об авторе}
\end{center}
\begin{enumerate}

    \item Рахматуллаев Музаффар Мухаммаджанович, канд. физ. мат. наук,
    старший научный сотрудник института математики
    при Национальном  Университете Узбекистана.
    \\адрес: Институт математики, ул. Дурмон йули, 29, Ташкент, 100125, Узбекистан.
    \\e-mail: mrahmatullaev@rambler.ru

\end{enumerate}

\begin{thebibliography}{99}


\bibitem{1} Rozikov U.A. \textit{Gibbs measures on Cayley trees}
(World scientific.,2013).
\bibitem{2} Блехер П.М.,Ганиходжаев Н.Н. \textit{Теория вероят. и ее
примен.},\textbf{35}(2), 920-930 (1990).
\bibitem{3} Zachary S. \textit{Ann. Prob.},\textbf{11}, 894-903 (1983).
\bibitem{4} Ганиходжаев Н.Н., Розиков У.А. \textit{ТМФ},\textbf{111}(1), 109-117 (1997).
\bibitem{GRS} Gandolfo, D., Ruiz, J., Shlosman, \textit{J.  Stat. Phys.} \textbf{153}(3), 400-411 (2013).
 \textit{ТМФ},\textbf{112}(1), 170-176 (1997).
\bibitem{5} Розиков У.А. \textit{ТМФ},\textbf{112}(1), 170-176 (1997).
\bibitem{6} Розиков У.А. \textit{ТМФ},\textbf{118}(1), 95-104 (1999).
\bibitem{7} Розиков У.А., Рахматуллаев М.М. \textit{ТМФ},\textbf{156}(2), 292-302
(2008).
\bibitem{8} Розиков У.А., Рахматуллаев М.М. \textit{ТМФ},\textbf{160}(3), 507-516
(2009).
\bibitem{9} Розиков У.А., Рахматуллаев М.М. \textit{Доклады АН РУз},\textbf{4},
12-15 (2008).
\bibitem{10} Розиков У.А., Рахматуллаев М.М. \textit{УзМЖ},\textbf{2},
144-152 (2009).
\bibitem{11} Ганиходжаев Н.Н. \textit{Доклады АН РУз},\textbf{4},
3-5 (1994).
\bibitem{*1} Я.Г. Синай,  {\it Теория фазовых переходов. Строгые
результаты}. (Наука.М., 1980).

\bibitem{*2} C.Preston, {\it Gibbs states on countable sets}. (Cambridge
Univ. Press,1974.)

\bibitem{*3} Малышев В.А., Минлос Р.А. \textit{Гиббсовские случайные поля} (Наука.М., 1985).

\bibitem{*4} Рюэль Д. \textit{Статистическая механика} (Мир.М., 1971).

\bibitem{*15} de Jongh L.J., Miedema A.R. \textit{Adv. Phys.},\textbf{23},
397-413 (1974).

\bibitem{*16} Фейман Р. \textit{Статистическая механика} (Мир.М., 1978).

\bibitem{*17} Kasteleyn P.W. \textit{Phase transitions. In: Fundamental Problems in Statistical Mechanics.II. Proc. of 2nd NUFFIC Int.Summer Course(Noordwijk, The Netherlands, 1967)} (New York, 1968).

\bibitem{*18} Zimm B., Doty P., Iso K. \textit{Proc. Nat. Acad.Sci.USA},\textbf{45},
160 (1959).

\bibitem{*19} Хуанг К. \textit{Статистическая механика} (Мир.М., 1966).




\end{thebibliography}
\end{document}